\documentclass[11pt]{article}
\usepackage{amssymb}
\usepackage{amsmath}
\usepackage[all]{xy}

\def\Mod{\mathop{\rm Mod}\nolimits}
\def\id{\mathop{\rm id}\nolimits}
\def\fd{\mathop{\rm fd}\nolimits}
\def\Ext{\mathop{\rm Ext}\nolimits}
\def\Tor{\mathop{\rm Tor}\nolimits}
\def\Hom{\mathop{\rm Hom}\nolimits}
\def\Im{\mathop{\rm Im}\nolimits}

\title{\Large \bf Injective Envelopes and (Gorenstein) Flat Covers}
\author{Edgar E. Enochs$^{1}$ and Zhaoyong Huang$^{2}$ \\
{\footnotesize  $^{1}$Department of Mathematics, University of
Kentucky, Lexington, KY 40506, USA;} \\ {\footnotesize
$^{2}$Department of Mathematics, Nanjing University, Nanjing 210093,
Jiangsu Province, P. R. China}
\\{\footnotesize (email: enochs@ms.uky.edu; huangzy@nju.edu.cn)}}
\date{}
\begin{document}
\baselineskip=18pt \maketitle

\begin{abstract}
We characterize left Noetherian rings in terms of the duality property of injective
preenvelopes and flat precovers. For a left and right Noetherian ring $R$, we prove
that the flat dimension of the injective envelope of any (Gorenstein) flat left
$R$-module is at most the flat dimension of the injective envelope of $_RR$. Then
we get that the injective envelope of $_RR$ is (Gorenstein) flat if and only if the
injective envelope of every Gorenstein flat left $R$-module is (Gorenstein) flat,
if and only if the injective envelope of every flat left $R$-module is (Gorenstein) flat,
if and only if the (Gorenstein) flat cover of every injective left $R$-module is injective,
and if and only if the opposite version of one of these conditions is satisfied.
\end{abstract}

\vspace{0.5cm}

{\bf 1. Introduction}

\vspace{0.2cm}

Throughout this paper, all rings are associative with identity. For
a ring $R$, we use $\Mod R$ (resp. $\Mod R^{op}$) to denote the
category of left (resp. right) $R$-modules.

It is well known that (pre)envelopes and (pre)covers of modules are
dual notions. These are fundamental and important in relative
homological algebra. Also note that coherent (resp. Noetherian)
rings can be characterized by the equivalence of the absolutely
purity (resp. injectivity) of modules and the flatness of their
character modules (see [ChS]). In this paper, under some conditions,
we show that we get a precover after applying a contravariant Hom
functor associated with a bimodule to a preenvelope. Then we obtain
an equivalent characterization of coherent (resp. Noetherian) rings
in terms of the duality property between absolutely pure (resp.
injective) preenvelopes and flat precovers.

Bican, El Bashir and Enochs proved in [BEE] that every module has a
flat cover for any ring. Furthermore, Enochs, Jenda and Lopez-Ramos
proved in [EJL] that every left module has a Gorenstein flat cover
for a right coherent ring. On the other hand, we know from [B] that
the injective envelope of $R$ is flat for a commutative Gorenstein
ring $R$. In [ChE] Cheatham and Enochs proved that for a commutative
Noetherian ring $R$, the injective envelope of $R$ is flat if and
only if the injective envelope of every flat $R$-module is flat. In
this paper, over a left and right Noetherian ring $R$, we will
characterize when the injective envelope of $_RR$ is (Gorenstein)
flat in terms of the (Gorenstein) flatness of the injective
envelopes of (Gorenstein) flat left $R$-modules and the injectivity
of the (Gorenstein) flat covers of injective left $R$-modules.

This paper is organized as follows.

In Section 2, we give some terminology and some preliminary results.
In particular, as generalizations of strongly cotorsion modules and
strongly torsionfree modules, we introduce the notions of
$n$-(Gorenstein) cotorsion modules and $n$-(Gorenstein) torsionfree
modules, and then establish a duality relation between right
$n$-(Gorenstein) torsionfree modules and left $n$-(Gorenstein)
cotorsion modules.

Let $R$ and $S$ be rings and let $_SU_R$ be a given $(S, R)$-bimodule.
For a subcategory $\mathcal{X}$ of $\Mod S$ (or $\Mod
R^{op}$), we denote by $\mathcal{X}^*=\{X^*\ |\ X\in \mathcal{X}\}$,
where $(-)^*=\Hom(-, {_SU_R})$. In Section 3, we prove that if
$\mathcal{C}$ is a subcategory of $\Mod S$ and $\mathcal{D}$ is a
subcategory of $\Mod R^{op}$ such that $\mathcal{C}^*\subseteq
\mathcal{D}$ and $\mathcal{D}^*\subseteq \mathcal{C}$, then a
homomorphism $f: A\to C$ in $\Mod S$ being a
$\mathcal{C}$-preenvelope of $A$ implies that $f^*: C^*\to A^*$ is a
$\mathcal{D}$-precover of $A^*$ in $\Mod R^{op}$. As applications of
this result, we get that a ring $R$ is left coherent if and only if
a monomorphism $f: A \rightarrowtail E$ in $\Mod R$ being an
absolutely pure (resp. injective) preenvelope of $A$ implies that $f^+: E^+
\twoheadrightarrow A^+$ is a flat precover of $A^+$ in $\Mod
R^{op}$, and that $R$ is left Noetherian if and only if a
monomorphism $f: A \rightarrowtail E$ in $\Mod R$ being an injective
preenvelope of $A$ is equivalent to $f^+: E^+ \twoheadrightarrow
A^+$ being a flat precover of $A^+$ in $\Mod R^{op}$, where $(-)^+$
is the character functor.

By using the result about the equivalent characterization of
Noetherian rings obtained in Section 3, we first prove in Section 4
that the flat dimension of the injective envelope of any
(Gorenstein) flat left $R$-module is at most the flat dimension of
the injective envelope of $_RR$ for a left and right Noetherian ring
$R$. Then we investigate the relation between the injective
envelopes of (Gorenstein) flat modules and (Gorenstein) flat covers
of injective modules. We give a list of equivalent conditions
relating these notions. For a left and right Noetherian ring $R$, we
prove that the injective envelope of $_RR$ is (Gorenstein) flat if
and only if the injective envelope of every Gorenstein flat left
$R$-module is (Gorenstein) flat, if and only if the injective
envelope of every flat left $R$-module is (Gorenstein) flat, if and
only if the (Gorenstein) flat cover of every injective left
$R$-module is injective, and if and only if the opposite version of
one of these conditions is satisfied. We remark that these
equivalent conditions hold true for commutative Gorenstein rings by [B],
but any of these conditions and the condition that $R$ is Gorenstein
are independent in general.

\vspace{0.5cm}

{\bf  2. Preliminaries}

\vspace{0.2cm}

In this section, we give some terminology and some preliminary
results for later use.

\vspace{0.2cm}

{\bf Definition 2.1.} ([E]) Let $R$ be a ring and $\mathcal{C}$ a
subcategory of $\Mod R$. The homomorphism $f: C\to D$ in $\Mod R$
with $C\in \mathcal{C}$ is said to be a {\it $\mathcal{C}$-precover}
of $D$ if for any homomorphism $g: C' \to D$ in $\Mod R$ with $C'\in
\mathcal{C}$, there exists a homomorphism $h: C'\to C$ such that the
following diagram commutes:
$$\xymatrix{ & C' \ar[d]^{g} \ar@{-->}[ld]_{h}\\
C \ar[r]^{f} & D}$$ The homomorphism $f: C\to D$ is said to be {\it
right minimal} if an endomorphism $h: C\to C$ is an automorphism
whenever $f=fh$. A $\mathcal{C}$-precover $f: C\to D$ is called a
{\it $\mathcal{C}$-cover} if $f$ is right minimal. Dually, the
notions of a {\it $\mathcal{C}$-preenvelope}, a {\it left minimal
homomorphism} and {\it a $\mathcal{C}$-envelope} are defined.

\vspace{0.2cm}

We begin with the following easy observation.

\vspace{0.2cm}

{\bf Lemma 2.2.} {\it Let $R$ be a ring and $\mathcal{D}$ a
subcategory of $\Mod R^{op}$, which is closed under direct products.
If $f_i: D_i \to M_i$ is a $\mathcal{D}$-precover of $M_i$ in $\Mod
R^{op}$ for any $i \in I$, where $I$ is an index set, then $\prod
_{i \in I}f_i: \prod _{i \in I}D_i\to \prod _{i \in I}M_i$ is a
$\mathcal{D}$-precover of $\prod _{i \in I}M_i$.}

\vspace{0.2cm}

{\it Proof.} It follows easily from the fact that $\Hom_R(D,-)$ commutes 
with direct products for any $D\in \mathcal{D}$. \hfill{$\square$}

\vspace{0.2cm}

As generalizations of flat modules and the flat dimension of
modules, the notions of Gorenstein flat modules and the Gorenstein
flat dimension of modules were introduced by Enochs, Jenda and
Torrecillas in [EJT] and by Holm in [H], respectively.

\vspace{0.2cm}

{\bf Definition 2.3.} ([EJT] and [H]) Let $R$ be a ring. A module $M$
in $\Mod R$ is called {\it Gorenstein flat} if there exists an exact
sequence:
$$\mathbb{F}: \cdots \to F_1 \to F_0 \to F^0 \to F^1 \to \cdots$$
in $\Mod R$ with all terms flat, such that $M=\Im (F_0 \to F^0)$ and
the sequence $I\bigotimes_R\mathbb{F}$ is exact for any injective right
$R$-module $I$. We use $\mathcal{GF}(R)$ to denote the subcategory
of $\Mod R$ consisting of Gorenstein flat modules. The {\it Gorenstein
flat dimension} of $M$ is defined as $\inf \{ n\ |\ $there exists an
exact sequence $0 \to G_n \to \cdots \to G_1 \to G_0 \to M \to 0$ in
$\Mod R$ with $G_i$ Gorenstein flat for any $0 \leq i \leq n\}$.

\vspace{0.2cm}

{\bf Lemma 2.4.} {\it For a ring $R$, an injective left (or right)
$R$-module is flat if and only if it is Gorenstein flat.}

\vspace{0.2cm}

{\it Proof.} The necessity is trivial. Notice that any Gorenstein
flat module can be embedded into a flat module, so it is easy to see
that a Gorenstein flat module is flat if it is injective.
\hfill{$\square$}

\vspace{0.2cm}

{\bf Lemma 2.5.} {\it ([H, Theorem 3.7]) For a right coherent ring
$R$, $\mathcal{GF}(R)$ is closed under extensions and under direct
summands.}

\vspace{0.2cm}

{\bf Definition 2.6.} ([X]) Let $R$ be a ring. A module $M$ in $\Mod
R$ is called {\it strongly cotorsion} if $\Ext_R^1(X, M)=0$ for any
$X\in \Mod R$ with finite flat dimension. A module $N$ in $\Mod
R^{op}$ is called {\it strongly torsionfree} if $\Tor_1^R(N, X)=0$
for any $X\in \Mod R$ with finite flat dimension.

\vspace{0.2cm}

We generalize these notions and introduce the notions of
$n$-cotorsion modules and $n$-torsionfree modules and that of
$n$-Gorenstein cotorsion modules and $n$-Gorenstein torsionfree
modules as follows.

\vspace{0.2cm}

{\bf Definition 2.7.} Let $R$ be a ring and $n$ a positive integer.

(1) A module $M$ in $\Mod R$ is called {\it $n$-cotorsion} if
$\Ext_R^1(X, M)=0$ for any $X\in \Mod R$ with flat dimension at most
$n$. A module $N$ in $\Mod R^{op}$ is called {\it $n$-torsionfree}
if $\Tor_1^R(N, X)=0$ for any $X\in \Mod R$ with flat dimension at
most $n$.

(2) A module $M$ in $\Mod R$ is called {\it $n$-Gorenstein
cotorsion} if $\Ext_R^1(X, M)=0$ for any $X\in \Mod R$ with
Gorenstein flat dimension at most $n$; and $M$ is called {\it
strongly Gorenstein cotorsion} if it is $n$-Gorenstein cotorsion for
all $n$. A module $N$ in $\Mod R^{op}$ is called {\it $n$-Gorenstein
torsionfree} if $\Tor_1^R(N, X)=0$ for any $X\in \Mod R$ with
Gorenstein flat dimension at most $n$; and $N$ is called {\it
strongly Gorenstein torsionfree} if it is $n$-Gorenstein torsionfree
for all $n$.

\vspace{0.2cm}

{\it Remark.} (1) We have the descending chains: $\{$1-cotorsion
modules$\}\supseteq \{$2-cotorsion modules$\}\supseteq \cdots
\supseteq \{$strongly cotorsion modules$\}$, and $\{$1-torsionfree
modules$\}\supseteq \{$2-torsionfree modules$\}\supseteq \cdots
\supseteq \{$strongly torsionfree modules$\}$. In particular,
$\{$strongly cotorsion modules$\}$\linebreak $=\bigcap _{n\geq
1}\{n$-cotorsion modules$\}$ and $\{$strongly torsionfree
modules$\}=\bigcap _{n\geq 1}\{n$-torsionfree modules$\}$.

(2) Similarly, we have the descending chains: $\{$1-Gorenstein
cotorsion modules$\}\supseteq \{$2-Gorenstein cotorsion
modules$\}\supseteq \cdots \supseteq \{$strongly Gorenstein
cotorsion modules$\}$, and $\{$1-Gorenstein torsionfree
modules$\}\supseteq \{$2-Gorenstein torsionfree modules$\}\supseteq
\cdots \supseteq \{$strongly Gorenstein torsionfree modules$\}$. In
particular, $\{$strongly Gorenstein cotorsion modules$\}=\bigcap
_{n\geq 1}\{n$-Gorenstein cotorsion modules$\}$ and $\{$strongly
Gorenstein torsionfree modules$\}=\bigcap _{n\geq 1}\{n$-Gorenstein
torsionfree modules$\}$.

\vspace{0.2cm}

We denote by $(-)^+=\Hom _{\mathbb{Z}}(-, \mathbb{Q}/\mathbb{Z})$,
where $\mathbb{Z}$ is the ring of integers and
$\mathbb{Q}$ is the ring of rational numbers. The
following result establishes a duality relation between right
$n$-(Gorenstein) torsionfree modules and left $n$-(Gorenstein)
cotorsion modules.

\vspace{0.2cm}

{\bf Proposition 2.8.} {\it Let $R$ be a ring and $N$ a module in
$\Mod R^{op}$. Then for any $n\geq 1$, we have

(1) $N$ is $n$-torsionfree if and only if $N^+$ is $n$-cotorsion. In
particular, $N$ is strongly torsionfree if and only if $N^+$ is
strongly cotorsion.

(2) $N$ is $n$-Gorenstein torsionfree if and only if $N^+$ is
$n$-Gorenstein cotorsion. In particular, $N$ is strongly Gorenstein
torsionfree if and only if $N^+$ is strongly Gorenstein cotorsion.}

\vspace{0.2cm}

{\it Proof.} By [CE, p.120, Proposition 5.1], we have that
$$[\Tor _1^R(N, A)]^+\cong \Ext _R^1(A, N^+)$$ for any $A\in \Mod
R$ and $N\in \Mod R^{op}$. Then both assertions follow easily.
\hfill{$\square$}

\vspace{0.2cm}

For a ring $R$, recall from [M] that a module $M$ in $\Mod R$ is
called {\it absolutely pure} if it is a pure submodule in every
module in $\Mod R$ that contains it, or equivalently, if it is pure
in every injective module in $\Mod R$ that contains it. Absolutely
pure modules are also known as {\it FP-injective modules}. It is
trivial that an injective module is absolutely pure. By [M, Theorem
3], a ring $R$ is left Noetherian if and only if every absolutely
pure module in $\Mod R$ is injective.

\vspace{0.2cm}

{\bf Lemma 2.9.} {\it (1) ([ChS, Theorem 1]) A ring $R$ is left
coherent if and only if a module $A$ in $\Mod R$ being absolutely
pure is equivalent to $A^+$ being flat in $\Mod R^{op}$.

(2) ([ChS, Theorem 2] and [F, Theorem 2.2])  A ring $R$ is left
(resp. right) Noetherian if and only if a module $E$ in $\Mod R$
(resp. $\Mod R^{op}$) being injective is equivalent to $E^+$ being
flat in $\Mod R^{op}$ (resp. $\Mod R$), and if and only if the
injective dimension of $M$ and the flat dimension of $M^+$ are
identical for any $M\in \Mod R$ (resp. $\Mod R^{op}$).}

\vspace{0.2cm}

As a generalization of injective modules, the notion of Gorenstein
injective modules was introduced by Enochs and Jenda in [EJ1] as
follows.

\vspace{0.2cm}

{\bf Definition 2.10.} ([EJ1]) Let $R$ be a ring. A module $M$ in
$\Mod R$ is called {\it Gorenstein injective} if there exists an
exact sequence:
$$\mathbb{I}:\cdots \to I_1 \to I_0 \to I^0 \to I^1 \to \cdots$$
in $\Mod R$ with all terms injective, such that $M=\Im (I_0 \to
I^0)$ and the sequence $\Hom _R(I,\mathbb{I})$ for any injective
left $R$-module $I$. We use $\mathcal{GI}(R)$ to denote the subcategory
of $\Mod R$ consisting of Gorenstein injective modules.

\vspace{0.5cm}

{\bf  3. The duality between preenvelopes and precovers}

\vspace {0.2cm}

In this section, we study the duality properties between
preenvelopes and precovers.

Let $R$ and $S$ be rings and let $_SU_R$ be a given $(S,R)$-bimodule.
We denote by $(-)^*=\Hom(-, {_SU_R})$. For a
subcategory $\mathcal{X}$ of $\Mod S$ (or $\Mod R^{op}$), we denote
by $\mathcal{X}^*=\{X^*\ |\ X\in \mathcal{X}\}$. For any $X\in \Mod
S$ (or $\Mod R^{op}$), $\sigma _X: X \to X^{**}$ defined by $\sigma
_X(x)(f)=f(x)$ for any $x\in X$ and $f\in X^*$ is the canonical
evaluation homomorphism.

\vspace {0.2cm}

{\bf Theorem 3.1.} {\it Let $\mathcal{C}$ be a subcategory of $\Mod
S$ and $\mathcal{D}$ a subcategory of $\Mod R^{op}$ such that
$\mathcal{C}^*\subseteq \mathcal{D}$ and $\mathcal{D}^*\subseteq
\mathcal{C}$. If $f: A\to C$ is a $\mathcal{C}$-preenvelope of a
module $A$ in $\Mod S$, then $f^*: C^*\to A^*$ is a
$\mathcal{D}$-precover of $A^*$ in $\Mod R^{op}$.}

\vspace{0.2cm}

{\it Proof.} Assume that $f: A\to C$ is a $\mathcal{C}$-preenvelope
of a module $A$ in $\Mod S$. Then we have a homomorphism $f^*:
C^*\to A^*$ in $\Mod R^{op}$ with $C^*\in \mathcal{C}^*\subseteq
\mathcal{D}$. Let $g: D\to A^*$ be a homomorphism in $\Mod R^{op}$
with $D \in \mathcal{D}$. Then $D^*\in \mathcal{D}^*\subseteq
\mathcal{C}$.

Consider the following diagram:
$$\xymatrix{A \ar[r]^f \ar[d]_{\sigma _A}& C\ar@{-->}[ldd]^{h}\\
A^{**} \ar[d]_{g^*}& \\
D^*& }$$ Because $f: A\to C$ is a $\mathcal{C}$-preenvelope of $A$,
there exists a homomorphism $h: C \to D^*$ such that the above
diagram commutes, that is, $hf=g^*\sigma _A$. Then we have $\sigma
_A^*g^{**}=f^*h^*$. On the other hand, we have the following
commutative diagram:
$$\xymatrix{D \ar[r]^g \ar[d]_{\sigma _D}& A^*\ar[d]^{\sigma _{A^{*}}}\\
D^{**} \ar[r]^{g^{**}}& A^{***}}$$ that is, we have $\sigma
_{A^{*}}g=g^{**}\sigma _D$. By [AF, Proposition 20.14], $\sigma
_A^*\sigma _{A^{*}}=1_{A^*}$. So we have that $g=1_{A^*}g=\sigma
_A^*\sigma _{A^{*}}g=\sigma _A^*g^{**}\sigma _D=f^*(h^*\sigma _D)$,
that is, we get a homomorphism $h^*\sigma _D: D \to C^*$ such that
the following diagram commutes:
$$\xymatrix{ & D \ar[d]^{g} \ar@{-->}[ld]_{h^*\sigma _D}\\
C^* \ar[r]^{f^*} & A^*}$$ Thus $f^*: C^*\to A^*$ is a
$\mathcal{D}$-precover of $A^*$. \hfill{$\square$}

\vspace{0.2cm}

In the rest of this section, we will give some applications of
Theorem 3.1.

For a ring $R$ and a subcategory $\mathcal{X}$ of $\Mod R$ (or $\Mod
R^{op}$), we denote by $\mathcal{X}^+=\{X^+\ |\ X\in \mathcal{X}\}$.

\vspace {0.2cm}

{\bf Corollary 3.2.} {\it (1) Let $\mathcal{C}$ be a subcategory of
$\Mod R$ and $\mathcal{D}$ a subcategory of $\Mod R^{op}$ such that
$\mathcal{C}^+\subseteq \mathcal{D}$ and $\mathcal{D}^+\subseteq
\mathcal{C}$. If $f: A\to C$ is a $\mathcal{C}$-preenvelope of a
module $A$ in $\Mod R$, then $f^+: C^+\to A^+$ is a
$\mathcal{D}$-precover of $A^+$ in $\Mod R^{op}$.

(2) Let $f: M\to N$ be a homomorphism in $\Mod R$. If $f^+: N^+\to
M^+$ is left (resp. right) minimal, then $f$ is right (resp. left)
minimal.}

\vspace{0.2cm}

{\it Proof.} (1) Notice that $\Hom_R(-, R^+)\cong (-)^+$ by the
adjoint isomorphism theorem, so the assertion is an immediate
consequence of Theorem 3.1.

(2) Assume that $f^+: N^+\to M^+$ is left minimal. If $f: M\to N$ is
not right minimal, then there exists an endomorphism $h: M \to M$
such that $f=fh$ but $h$ is not an automorphism. So $f^+=h^+f^+$ and
$h^+$ is not an automorphism. It follows that $f^+$ is not left
minimal, which is a contradiction. Thus we conclude that $f$ is
right minimal. Similarly, we get that $f$ is left minimal if $f^+$
is right minimal.  \hfill{$\square$}

\vspace{0.2cm}

We use $\mathcal{AP}(R)$, $\mathcal{F}(R)$ and $\mathcal{I}(R)$ to
denote the subcategories of $\Mod R$ consisting of absolutely pure
modules, flat modules and injective modules, respectively. Recall
that an $\mathcal{F}(R)$-preenvelope is called a {\it flat
preenvelope}. By [E, Proposition 5.1] or [EJ2, Proposition 6.5.1],
we have that $R$ is a left coherent ring if and only if every module
in $\Mod R^{op}$ has a flat preenvelope. Also recall that an
$\mathcal{AP}(R)$-precover and an $\mathcal{I}(R)$-precover are
called an {\it absolutely pure precover} and an {\it injective
precover}, respectively. It is known that every module in $\Mod R$
has an absolutely pure precover for a left coherent ring $R$ by [P,
Theorem 2.6], and every module in $\Mod R$ has an injective precover
for a left Noetherian ring $R$ by [E, Proposition 2.2].

\vspace{0.2cm}

{\bf Corollary 3.3.} {\it Let $R$ be a left coherent ring. If a
homomorphism $f: A \to F$ in $\Mod R^{op}$ is a flat preenvelope of
$A$, then $f^+: F^+ \to A^+$ is an absolutely pure precover and an
injective precover of $A^+$ in $\Mod R$.}

\vspace{0.2cm}

{\it Proof.} By [F, Theorem 2.1] and Lemma 2.9(1) and the opposite version of
Corollary 3.2(1). \hfill{$\square$}

\vspace{0.2cm}

Recall that a $\mathcal{GI}(R)$-precover (resp. preenvelope) and a
$\mathcal{GF}(R)$-(pre)cover (resp. preenvelope) are called a {\it
Gorenstein injective precover} (resp. {\it preenvelope}) and a {\it
Gorenstein flat (pre)cover} (resp. {\it preenvelope}), respectively.
Let $R$ be a Gorenstein ring. Then every module in $\Mod R$ has a
Gorenstein injective preenvelope and a Gorenstein flat preenvelope
by [EJ2, Theorems 11.2.1 and 11.8.2], and every module in $\Mod
R^{op}$ has a Gorenstein injective precover and a Gorenstein flat
precover by [EJ2, Theorems 11.1.1] and [EJL, Theorem 2.11],
respectively.

\vspace{0.2cm}

{\bf Corollary 3.4.} {\it Let $R$ be a Gorenstein ring.

(1) If a monomorphism $f: A \rightarrowtail Q$ in $\Mod R$ is a
Gorenstein injective preenvelope of $A$, then $f^+: Q^+
\twoheadrightarrow A^+$ is a Gorenstein flat precover of $A^+$ in
$\Mod R^{op}$.

(2) If a homomorphism $f: A \to G$ in $\Mod R$ is a Gorenstein flat
preenvelope of $A$, then $f^+: G^+ \to A^+$ is a Gorenstein
injective precover of $A^+$ in $\Mod R^{op}$.}

\vspace{0.2cm}

{\it Proof.} By [H, Theorem 3.6], [EJ2, Corollary 10.3.9] and Corollary 3.2(1).
\hfill{$\square$}

\vspace{0.2cm}

Recall that an $\mathcal{AP}(R)$-preenvelope is called an {\it
absolutely pure preenvelope}. By [EJ2, Proposition 6.2.4], every
module in $\Mod R$ has an absolutely pure preenvelope. Recall that
an $\mathcal{I}(R)$-(pre)envelope and an $\mathcal{F}(R)$-(pre)cover
are called an {\it injective (pre)envelope} and a {\it flat
(pre)cover}, respectively. By [BEE, Theorem 3], every module in
$\Mod R$ has a flat cover. In the following result, we give an
equivalent characterization of left coherent rings in terms of the
duality property between absolutely pure preenvelopes and flat
precovers.

\vspace{0.2cm}

{\bf Theorem 3.5.} {\it The following statements are equivalent for a ring
$R$.

(1) $R$ is a left coherent ring.

(2) If a monomorphism $f: A \rightarrowtail E$ in $\Mod R$ is an
absolutely pure preenvelope of $A$, then $f^+: E^+
\twoheadrightarrow A^+$ is a flat precover of $A^+$ in $\Mod
R^{op}$.

(3) If a monomorphism $f: A \rightarrowtail E$ in $\Mod R$ is an
injective preenvelope of $A$, then $f^+: E^+ \twoheadrightarrow A^+$
is a flat precover of $A^+$ in $\Mod R^{op}$.}

\vspace{0.2cm}

{\it Proof.} $(1)\Rightarrow (2)$ and $(1)\Rightarrow (3)$ follow
from [F, Theorem 2.1], Lemma 2.9(1) and Corollary 3.2(1).

$(2)\Rightarrow (1)$ Note that a left $R$-module $E$ is absolutely
pure if $E^+$ is flat (see [ChS]). So by (2), we have that a left
$R$-module $E$ is absolutely pure if and only if $E^+$ is flat. Then
it follows from Lemma 2.9(1) that $R$ is a left coherent ring.

$(3)\Rightarrow (2)$ Let $E\in \Mod R$ be absolutely pure. Then the inclusion
$E\hookrightarrow E^0(E)$ is pure exact, and so $[E^0(E)]^+\twoheadrightarrow E^+$
admits a section by [EJ2, Proposition 5.3.8]. It yields that $E^+$ is isomorphic
to a direct summand of $[E^0(E)]^+$. Note that $[E^0(E)]^+$ is flat by (3). So
$E^+$ is also flat. On the other hand, the character module of any flat module
in $\Mod R^{op}$ is injective (and hence absolutely pure) by [F, Theorem 2.1].
Thus we get the assertion by Corollary 3.2(1). \hfill{$\square$}

\vspace{0.2cm}

The following example illustrates that for a left (and right)
coherent ring $R$, neither of the converses of (2) and (3) in
Theorem 3.5 hold true in general.

\vspace{0.2cm}

{\bf Example 3.6.} Let $R$ be a von Neumann regular ring but not a
semisimple Artinian ring. Then $R$ is a left and right coherent ring
and every left $R$-module is absolutely pure by [M, Theorem 5]. So
both (2) and (3) in Theorem 3.5 hold true. Because $R$ is not a
semisimple Artinian ring, there exists a non-injective left
$R$-module $M$. Then we have

(1) The injective envelope $\alpha: M \hookrightarrow E^0(M)$ is not
an absolutely pure preenvelope, because there does not exist a
homomorphism $E^0(M) \to M$ such that the following diagram
commutes:
$$\xymatrix{M \ar[r]^{\alpha} \ar[d]_{1_M}& E^0(M)\\
M & }$$ So the converse of (2) in Theorem 3.5 does not hold true.

(2) Note that the identity homomorphism $1_M: M \to M$ is an
absolutely pure envelope of $M$ but not an injective preenvelope of
$M$. So the converse of (3) in Theorem 3.5 does not hold true.

Thus we also get that neither of the converses of Corollary 3.2(1)
and Theorem 3.1 hold true in general.

\vspace{0.2cm}

In the following result, we give an equivalent characterization of
left Noetherian rings in terms of the duality property between
injective preenvelopes and flat precovers.

\vspace{0.2cm}

{\bf Theorem 3.7.} {\it The following statements are equivalent for a ring
$R$.

(1) $R$ is a left Noetherian ring.

(2) A monomorphism $f: A \rightarrowtail E$ in $\Mod R$ is an
injective preenvelope of $A$ if and only if $f^+: E^+
\twoheadrightarrow A^+$ is a flat precover of $A^+$ in $\Mod
R^{op}$.

(3) $R$ is a left coherent ring, and a monomorphism $f: A \rightarrowtail E$
is an injective envelope of $A$ if $f^+: E^+ \twoheadrightarrow A^+$ is a
flat cover of $A^+$.}

\vspace{0.2cm}

{\it Proof.} $(2)\Rightarrow (1)$ follows from Lemma 2.9(2).

$(1)\Rightarrow (2)$ Assume that $R$ is a left Noetherian ring. If
$f: A \rightarrowtail E$ in $\Mod R$ is an injective preenvelope of
$A$, then $f^+: E^+ \twoheadrightarrow A^+$ is a flat precover of
$A^+$ by Theorem 3.5. Conversely, if $f^+: E^+ \twoheadrightarrow
A^+$ is a flat precover of $A^+$ in $\Mod R^{op}$, then $E$ is
injective by Lemma 2.9(2), and so $f: A \rightarrowtail E$ is an
injective preenvelope of $A$.

$(2)\Rightarrow (3)$ follows from Theorem 3.5 and Corollary 3.2(2).

$(3)\Rightarrow (1)$ By (3) and Lemma 2.9(1), we have that a module $E$ in
$\Mod R$ is injective if and only if $E^+$ is flat in $\Mod R^{op}$. So
from Lemma 2.9(2) we get the assertion. \hfill{$\square$}

\vspace {0.5cm}

{\bf 4. Injective envelopes of (Gorenstein) flat modules}

\vspace{0.2cm}

In this section, $R$ is a left and right Noetherian ring. We will
investigate the relation between the injective envelopes of
(Gorenstein) flat modules and (Gorenstein) flat covers of injective
modules.

For a module $M$ in $\Mod R$, we denote the injective envelope and
the flat cover of $M$ by $E^0(M)$ and $F_0(M)$ respectively, and
denote the injective dimension and the flat dimension of $M$ by $\id
_RM$ and $\fd _RM$ respectively.

\vspace{0.2cm}

{\bf Theorem 4.1.} \begin{align*}
&\ \ \ \ \fd _RE^0(_RR)\\
& =\sup\{\fd _RE^0(F)\ |\ F\in\Mod R \ \text{\it is flat}\}\\
& =\sup\{\fd _RE^0(G)\ |\ G\in \Mod R \ \text{\it is Gorenstein flat}\}.
\end{align*}

\vspace{0.2cm}

{\it Proof.} It is trivial that $\fd _RE^0(_RR)\leq \sup\{\fd
_RE^0(F)\ |\ F\in\Mod R$ is flat$\}\leq$\linebreak $\sup\{\fd
_RE^0(G)\ |\ G\in \Mod R$ is Gorenstein flat$\}$. So it suffices to
prove the opposite inequalities.

We first prove $\fd _RE^0(_RR)\geq \sup\{\fd _RE^0(F)\ |\ F\in\Mod
R$ is flat$\}$. Without loss of generality, suppose $\fd
_RE^0(_RR)=n<\infty$. Because $R$ is a left Noetherian ring and
$_RR\hookrightarrow E^0(_RR)$ is an injective envelope of $_RR$,
$[E^0(_RR)]^+\twoheadrightarrow (_RR)^+$ is a flat precover of
$(_RR)^+$ in $\Mod R^{op}$ by Theorem 3.7.

Let $F\in \Mod R$ be flat. Then $F^+$ is injective in $\Mod R^{op}$
by [F, Theorem 2.1]. Because $(_RR)^+$ is an injective cogenerator for
$\Mod R^{op}$ by [S, p.32, Proposition 9.3], $F^+$ is isomorphic to
a direct summand of $\prod _{i\in I}(_RR)^+$ for some set $I$. So
$F_0(F^+)$ is isomorphic to a direct summand of any flat precover of
$\prod _{i\in I}(_RR)^+$. Notice that $\mathcal{F}(R^{op})$ is
closed under direct products by [C, Theorem 2.1], so $\prod _{i\in
I}[E^0(_RR)]^+ \twoheadrightarrow \prod _{i\in I}(_RR)^+$ is a flat
precover of $\prod _{i\in I}(_RR)^+$ by the above argument and Lemma
2.2. Thus we get that $F_0(F^+)$ is isomorphic to a direct summand
of $\prod _{i\in I}[E^0(_RR)]^+$. Because $\fd _RE^0(_RR)=n$, $\id
_{R^{op}}[E^0(_RR)]^+=n$ by [F, Theorem 2.1]. So $\id _{R^{op}}\prod
_{i\in I}[E^0(_RR)]^+=n$ and $\id _{R^{op}}F_0(F^+)\leq n$.

Because $R$ is a right Noetherian ring, we get an injective
preenvelope $F^{++}\rightarrowtail [F_0(F^+)]^+$ of $F^{++}$ with
$\fd _R[F_0(F^+)]^+\leq n$ by Lemma 2.9(2). By [S, p.48, Exercise
41], there exists a monomorphism $F \rightarrowtail F^{++}$. So
$E^0(F)$ is isomorphic to a direct summand of $E^0(F^{++})$, and
hence a direct summand of $[F_0(F^+)]^+$. It follows that $\fd
_RE^0(F)\leq n$. Thus we get that $\fd _RE^0(_RR)\geq \sup\{\fd
_RE^0(F)\ |\ F\in\Mod R$ is flat$\}$.

Next, we prove $\sup\{\fd _RE^0(F)\ |\ F\in\Mod R$ is flat$\}\geq
\sup\{\fd _RE^0(G)\ |\ G\in \Mod R$ is Gorenstein flat$\}$. Let
$G\in \Mod R$ be Gorenstein flat. Then $G$ can be embedded into a
flat left $R$-module $F$, which implies that $E^0(G)$ is isomorphic
to a direct summand of $E^0(F)$ and $\fd _RE^0(G) \leq \fd
_RE^0(F)$. Thus the assertion follows. \hfill{$\square$}

\vspace{0.2cm}

{\bf Corollary 4.2.} {\it The following statements are equivalent.

(1) $E^0(_RR)$ is flat.

(2) $E^0(F)$ is flat for any flat left $R$-module $F$.

(3) $E^0(_RR)$ is Gorenstein flat.

$(i)^{op}$ The opposite version of (i) ($1 \leq i \leq 3$).}

\vspace{0.2cm}

{\it Proof.} $(2) \Leftrightarrow (1) \Leftrightarrow (3)$ follow
from Theorem 4.1 and Lemma 2.4, and $(1) \Rightarrow (1)^{op}$
follows from [Mo, Theorem 1]. Symmetrically, we get the opposite
versions of the implications mentioned above. \hfill{$\square$}

\vspace{0.2cm}

We know from [EJL, Theorem 2.12] that every module in $\Mod R$ has a
Gorenstein flat cover. For a module $M$ in $\Mod R$, we denote the
Gorenstein flat cover of $M$ by $GF_0(M)$.

\vspace{0.2cm}

{\bf Theorem 4.3.} {\it The following statements are equivalent.

(1) $E^0(_RR)$ is Gorenstein flat.

(2) $E^0(F)$ is Gorenstein flat for any flat left $R$-module $F$.

(3) $E^0(G)$ is Gorenstein flat for any Gorenstein flat left
$R$-module $G$.

(4) $GF_0(M)$ is injective for any 1-Gorenstein cotorsion left
$R$-module $M$.

(5) $GF_0(M)$ is injective for any strongly Gorenstein cotorsion
left $R$-module $M$.

(6) $GF_0(E)$ is injective for any injective left $R$-module $E$.

(7) $E^0(N)$ is flat for any 1-Gorenstein torsionfree right
$R$-module $N$.

(8) $E^0(N)$ is Gorenstein flat for any 1-Gorenstein torsionfree
right $R$-module $N$.

(9) $E^0(N)$ is flat for any strongly Gorenstein torsionfree right
$R$-module $N$.

(10) $E^0(N)$ is Gorenstein flat for any strongly Gorenstein
torsionfree right $R$-module $N$.}

\vspace{0.2cm}

{\it Proof.} $(1)\Rightarrow (2)$ By Corollary 4.2.

$(2)\Rightarrow (3)$ Let $G\in \Mod R$ be Gorenstein flat. From the
proof of Theorem 4.1, we know that there exists a flat left
$R$-module $F$ such that $E^0(G)$ is isomorphic to a direct summand
of $E^0(F)$. By (2), $E^0(F)$ is Gorenstein flat. Because
$\mathcal{GF}(R)$ is closed under direct summands by Lemma 2.5,
$E^0(G)$ is also Gorenstein flat.

$(3)\Rightarrow (4)$ Let $M\in \Mod R$ be 1-Gorenstein cotorsion.
Consider the following push-out diagram:
$$\xymatrix{& & 0 \ar[d] & 0 \ar[d] & \\
0 \ar[r] & X \ar[r] \ar@{=}[d]
& GF_0(M) \ar[r] \ar[d] & M \ar[r] \ar[d] & 0\\
0 \ar[r] & X \ar[r] & E^0(GF_0(M)) \ar[r] \ar[d] & N \ar[r] \ar[d] & 0\\
& & T \ar@{=}[r] \ar[d] & T \ar[d] &\\
& & 0 & 0 &}$$ By (3), $E^0(GF_0(M))$ is Gorenstein flat. Then by
the exactness of the middle column in the above diagram, we have
that the Gorenstein flat dimension of $T$ is at most one. So $\Ext
_R^1(T, M)=0$ and the rightmost column $0\to M \to N \to T \to 0$ in
the above diagram splits, which implies that $M$ is isomorphic to a
direct summand of $N$. It follows that $GF_0(M)$ is isomorphic to a
direct summand of $GF_0(N)$.

Note that $\mathcal{GF}(R)$ is closed under extensions by Lemma 2.5.
So $\Ext _R^1(G, X)=0$ for any Gorenstein flat left $R$-module $G$
by [X, Lemma 2.1.1], which implies that $E^0(GF_0(M))$ is a Gorenstein flat
precover of $N$ by the exactness of the middle row in the above
diagram. Thus $GF_0(N)$ is isomorphic to a direct summand of
$E^0(GF_0(M))$, and therefore $GF_0(N)$ and $GF_0(M)$ are injective.

$(4)\Rightarrow (5)\Rightarrow (6)$ are trivial.

$(6)\Rightarrow (1)$ Consider the following pull-back diagram:
$$\xymatrix{&  0 \ar[d] & 0 \ar[d] & &\\
&  W \ar@{=}[r] \ar[d]
& W \ar[d] & & \\
0 \ar[r] & Y \ar[r] \ar[d] & GF_0(E^0(_RR)) \ar[r] \ar[d] & H \ar[r] \ar@{=}[d] & 0\\
0 \ar[r] & _RR \ar[r] \ar[d] & E^0(_RR) \ar[r] \ar[d] & H \ar[r] & 0\\
& 0 & 0 & &}$$ Then $_RR$ is isomorphic to a direct summand of $Y$,
and so $E^0(_RR)$ is isomorphic to a direct summand of $E^0(Y)$.
Because $GF_0(E^0(_RR))$ is injective by (6), the exactness of the
middle row in the above diagram implies that $E^0(Y)$ is isomorphic
to a direct summand of $GF_0(E^0(_RR))$. Thus $E^0(_RR)$ is also
isomorphic to a direct summand of $GF_0(E^0(_RR))$. Again by Lemma
2.5, $\mathcal{GF}(R)$ is closed under direct summands, so we get
that $E^0(_RR)$ is Gorenstein flat.

$(1)+(4)\Rightarrow (7)$ Let $N\in \Mod R^{op}$ be 1-Gorenstein
torsionfree. Then $N^+\in \Mod R$ is 1-Gorenstein cotorsion by
Proposition 2.8, and so $GF_0(N^+)$ is injective by (4).

From the epimorphism $GF_0(N^+)\twoheadrightarrow N^+$ in $\Mod R$,
we get a monomorphism $ N^{++}\rightarrowtail [GF_0(N^+)]^+$ in
$\Mod R^{op}$ with $[GF_0(N^+)]^+$ flat by Lemma 2.9(2). Then
$E^0(N^{++})$ is isomorphic to a direct summand of
$E^0([GF_0(N^+)]^+)$. On the other hand, there exists a monomorphism
$N\rightarrowtail N^{++}$ by [S, p.48, Exercise 41], $E^0(N)$ is
isomorphic to a direct summand of $E^0(N^{++})$, and hence a direct
summand of $E^0([GF_0(N^+)]^+)$. Because $E^0([GF_0(N^+)]^+)$ is
flat by (1) and Corollary 4.2, $E^0(N)$ is also flat.

$(7) \Rightarrow (8) \Rightarrow (10)$ and $(7)\Rightarrow (9)
\Rightarrow (10)$ are trivial.

$(10) \Rightarrow (1)$ By (10), we have that $E^0(R_R)$ is
Gorenstein flat. Then $E^0(_RR)$ is also Gorenstein flat by
Corollary 4.2. \hfill{$\square$}

\vspace{0.2cm}

The following result is an analogue of Theorem 4.3.

\vspace{0.2cm}

{\bf Theorem 4.4.} {\it The following statements are equivalent.

(1) $E^0(_RR)$ is flat.

(2) $E^0(F)$ is flat for any flat left $R$-module $F$.

(3) $E^0(G)$ is flat for any Gorenstein flat left $R$-module $G$.

(4) $F_0(M)$ is injective for any 1-cotorsion left $R$-module $M$.

(5) $F_0(M)$ is injective for any strongly cotorsion left $R$-module
$M$.

(6) $F_0(E)$ is injective for any injective left $R$-module $E$.

(7) $E^0(N)$ is flat for any 1-torsionfree right $R$-module $N$.

(8) $E^0(N)$ is Gorenstein flat for any 1-torsionfree right
$R$-module $N$.

(9) $E^0(N)$ is flat for any strongly torsionfree right $R$-module
$N$.

(10) $E^0(N)$ is Gorenstein flat for any strongly torsionfree right
$R$-module $N$.}

\vspace{0.2cm}

{\it Proof.} The proof is similar to that of Theorem 4.3, so we omit
it. \hfill{$\square$}

\vspace{0.2cm}

Putting the results in this section and their opposite versions
together we have the following

\vspace{0.2cm}

{\bf Theorem 4.5.} {\it The following statements are equivalent.

(1) $E^0(_RR)$ is flat.

(2) $E^0(G)$ is flat for any (Gorenstein) flat left $R$-module $G$.

(3) $E^0(M)$ is (Gorenstein) flat for any 1-torsionfree left
$R$-module $M$.

(4) $E^0(M)$ is (Gorenstein) flat for any strongly torsionfree left
$R$-module $M$.

(5) $F_0(M)$ is injective for any 1-cotorsion left $R$-module $M$.

(6) $F_0(M)$ is injective for any strongly cotorsion left $R$-module
$M$.

(7) $F_0(E)$ is injective for any injective left $R$-module $E$.

$(i)^{op}$ The opposite version of (i) ($1 \leq i \leq 7$).

(G1) $E^0(_RR)$ is Gorenstein flat.

(G2) $E^0(G)$ is Gorenstein flat for any (Gorenstein) flat left
$R$-module $G$.

(G3) $E^0(M)$ is (Gorenstein) flat for any 1-Gorenstein torsionfree
left $R$-module $M$.

(G4) $E^0(M)$ is (Gorenstein) flat for any strongly Gorenstein
torsionfree left $R$-module $M$.

(G5) $GF_0(M)$ is injective for any 1-Gorenstein cotorsion left
$R$-module $M$.

(G6) $GF_0(M)$ is injective for any strongly Gorenstein cotorsion
left $R$-module $M$.

(G7) $GF_0(E)$ is injective for any injective left $R$-module $E$.

$(Gi)^{op}$ The opposite version of (Gi) ($1 \leq i \leq 7$).}

\vspace{0.2cm}

By [B, The Fundamental Theorem], $E^0(R)$ is flat if $R$ is a commutative
Gorenstein ring. In addition, there exists a non-commutative Gorenstein ring
$R$ such that $E^0(_RR)$ is flat (see [IM, Section 3]). However, the conditions
``$E^0(_RR)$ is flat" and ``$R$ is Gorenstein" are independent in general
as shown in the following example.

\vspace{0.2cm}

{\bf Example 4.6.} (1) Let $R$ be a finite-dimensional algebra given by the quiver:
$$\xymatrix{1 \ar@(ul,dl)_{\alpha}\ar[r]^{\beta} & 2 \ar @<2pt> [r]^{\gamma} &
3 \ar @<2pt> [l]^{\delta}}$$
modulo the ideal generated by $\{\alpha^2, \gamma\beta\alpha, \gamma\delta,
\beta\alpha-\delta\gamma\beta\}$. Then $E^0(_RR)$ is
flat and $R$ is not Gorenstein.

%In addition, let $R=K[x^3, x^4, x^5]$ with $K$ a field. Then $E^0(R)$ is
%flat, but $R$ is not Gorenstein.

(2) Let $R$ be a finite-dimensional algebra given by the quiver:
$$2 \longleftarrow 1 \longrightarrow 3.$$
Then $R$ is Gorenstein with $\id_RR=\id_{R^{op}}R=1$
and $\fd _RE^0(_RR)=1$.

\vspace{0.2cm}

Consider the following conditions for any $n\geq 0$.

(1) $\fd_RE^0(_RR)\leq n$.

(2) $\id_R F_0(E)\leq n$ for any injective left $R$-module $E$.

When $n=0$, $(1)\Leftrightarrow (2)$ by Theorem 4.5. However, when
$n\geq 1$, neither ``$(1)\Rightarrow (2)$" nor ``$(2)\Rightarrow
(1)$" hold true in general as shown in the following example.

\vspace{0.2cm}

{\bf Example 4.7.} Let $R$ be a finite-dimensional algebra over a
field $K$ and $\Delta$ the quiver:
$$\xymatrix{ 1 \ar @<2pt> [r]^{\alpha} &
2 \ar @<2pt> [l]^{\beta} \ar[r]^{\gamma} & 3.}
$$

(1) If $R=K\Delta /(\alpha \beta \alpha)$, then $\fd_RE^0(_RR)=1$
and $\fd_{R^{op}}E^0(R_R)\geq 2$. We have that
$\mathbb{D}[E^0(R_R)]\twoheadrightarrow \mathbb{D}(R_R)$ is the flat
cover of the injective left $R$-module $\mathbb{D}(R_R)$ with
$\id_R\mathbb{D}[E^0(R_R)]$ \linebreak $\geq 2$, where
$\mathbb{D}(-)=\Hom_K(-, K)$.

(2) If $R=K\Delta /(\gamma \alpha , \beta \alpha)$, then
$\fd_RE^0(_RR)=2$ and $\fd_{R^{op}}E^0(R_R)=1$. We have that
$\mathbb{D}[E^0(R_R)]\twoheadrightarrow \mathbb{D}(R_R)$ is the flat
cover of the injective left $R$-module $\mathbb{D}(R_R)$ with
$\id_R\mathbb{D}[E^0(R_R)]=1$. Because $\mathbb{D}(R_R)$ is an
injective cogenerator for $\Mod R$, $\id _RF_0(E)\leq 1$ for any
injective left $R$-module $E$.

\vspace{0.5cm}

{\bf Acknowledgements.} Part of this paper was written while the second
author was visiting University of Kentucky from June to August,
2009. The second author was partially supported by the
Specialized Research Fund for the Doctoral Program of Higher
Education (Grant No. 20100091110034) and NSF
of Jiangsu Province of China (Grant Nos. BK2010047, BK2010007).
The authors thank the referee for the useful suggestions.

\end{document}